\documentclass[a4paper,12pt]{amsart}

\usepackage{amssymb,amsbsy,amsmath,amsfonts,amssymb,amscd}
\usepackage{latexsym}
\usepackage{graphics}
\usepackage{color}
\input xy
\xyoption{all}

\def\eea{\end{eqnarray*}}
\def\bea{\begin{eqnarray*}}

\newcommand\dual{\mathrel{\raise3pt\hbox{$\underline{\mathrm{\thinspace d
\thinspace}}$}}}
\newcommand\qe{\ifhmode\unskip\nobreak\fi\quad $\Box$}       

\def\BOX{\hfill\lower.5\baselineskip\hbox{$\Box$}}

\newtheorem{theorem}{Theorem}[section]

\newtheorem{proposition}[theorem]{Proposition}

\newtheorem{theo}{Theorem}[section]
\newtheorem{remarkk}[theo]{Remark}

\newtheorem{defin}[theo]{Definition}

\newtheorem{example}[theo]{Example}

\begin{document}

\title[Corrigendum]{ Corrigendum for ``A geometric proof of the Karpelevich-Mostow Theorem''}
\author{ Antonio J. Di Scala - Carlos Olmos}

\thanks{The first author was partially supported by GNSAGA of INdAM, PRIN 07 ``Differential Geometry and Global Analysis'' and MIUR of Italy.
The second author was supported by FaMAF-Universidad Nacional de C\'ordoba, CIEM-Conicet, Secyt-UNC y ANCyT, Argentina}



\begin{abstract}
Corollary 2.3 in our paper \cite{nostro} is false.
Here we show how to avoid the use of this corollary to give a simpler proof of Karpelevich-Mostow theorem.
We also include a short discussion of the original proof by Karpelevich.
\end{abstract}

\maketitle

We want to point out that Corollary 2.3 in our paper \cite{nostro} is false.
A counterexample is as follows. Let $M=G/K$ be a riemannian symmetric space of rank $> 1$ and
let $K.v_p$, $v_p \neq 0$, be a focal orbit of the action of $K$ on $T_pM$, where $p = [K]$. The normal space $\mathfrak{r} := \nu_{v_p}(K.v_p)$ is a Lie triple system having a non-trivial flat factor $\mathfrak{a}$. So the group of isometries of the totally geodesic submanifold $Exp(\mathfrak{r})$ contains as a factor a semisimple subgroup $G$ which preserves all the points at infinity produced by the geodesics with initial directions at $\mathfrak{a}$.
So such a $G$ is a counterexample to Corollary 2.3. The confusion comes from the rank one case in which Corollary 2.3 is true.

Since Corollary 2.3 was used in the proof of Theorem 4.1 it follows that the geometric proof of Karpelevich-Mostow theorem in \cite{nostro} does not work. Fortunately, Theorem 4.1 in \cite{nostro} can be avoided and actually {\bf the main idea of our proof still works giving an even simpler proof of the Karpelevich-Mostow theorem.} The key observation is the following proposition.

\begin{proposition}\label{out-in} Let $\mathcal{P}$ be a symmetric space of non positive curvature and let $G \subset \mathrm{Iso}(\mathcal{P})$ be a connected Lie subgroup of isometries. Let $M \subset \mathcal{P}$ be a totally geodesic submanifold invariant by $G$, i.e. if $x \in M$ and $g \in G$ then $g.x \in M$.
If $G$ has a totally geodesic orbit in $\mathcal{P}$ then $G$ has also a totally geodesic orbit in $M$.
\end{proposition}

\it Proof. \rm For the sake of completeness we repeat here part of the proof of Proposition 3.1 of our article \cite{nostro}.
Let $M' = G.p'$ be a totally geodesic orbit of the action of $G$ on $\mathcal{P}$. We assume $ p' \notin M$ otherwise $G$ has a totally geodesic orbit in $M$ and we have nothing to prove.
Let now $\gamma(t)$ be a unit speed geodesic of $X$ starting at $p'=\gamma(0) \in M'$ which realizes the distance $d$ between $p'$ and $M$. Observe that $\gamma'(d)$ is perpendicular to $M$ so $\gamma'(d)$ is also perpendicular to the $G$-orbit $G.\gamma(d)$. By elementary properties of Killing vector fields it follows that $\gamma'(t)$ is always perpendicular to the orbits $G.\gamma(t)$.

{ \bf We claim that $N = G.\gamma(d)$ is a totally geodesic submanifold of $M$.}

For each Killing vector field $X \in \mathfrak{g}$
define
\[ f(t) := \langle \nabla_{\gamma'(t)} X , X \rangle = \frac{\mathrm{d}}{\mathrm{d}t} \frac{||X||^2}{2} \]
\noindent
where $\nabla$ is the Levi-Civita connection of $S$.\\

Notice that $f(0)=f(d)=0$. Indeed, $f(0)=0$ follows from the Killing equation and the fact that $G.\gamma(0)$ is a totally geodesic submanifold of $S$. The equality $f(d)=0$ follows since $\gamma'(d)$ is perpendicular to the totally geodesic submanifold $M$.
Then by computing  $f'(t)$ we get:
\[ \begin{aligned} \frac{d}{dt}f(t) &= \langle \nabla_{\gamma'(t)} \nabla_{\gamma'(t)} X , X \rangle + ||\nabla_{\gamma'(t)} X||^2  \\
&= -\langle R_{X, \gamma'(t)}\gamma'(t),X \rangle + ||\nabla_{\gamma'(t)} X||^2 \\
&\geq 0 \, \,.
\end{aligned} \]
\noindent
Here $R$ is the curvature tensor of $S$ and the second equality follows from the fact that
Killing vector fields are Jacobi fields along geodesics.

Then $\frac{d}{dt}f(t) \equiv 0$ and so we have

\begin{equation}\label{dos-pajaros} \nabla_{\gamma'(t)} X \equiv  \langle R_{X, \gamma'(t)}\gamma'(t),X \rangle \equiv 0 \, \, .
\end{equation}

Notice that the second equality implies $R_{X, \gamma'(t)}(\cdot) \equiv 0$ since $S$ is a symmetric space.

Let $\eta(t)$ be a $\nabla$-parallel vector field along $\gamma(t)$ such that $\eta(0)$ is perpendicular to $G.\gamma(0)$.
To complete the proof that $G.\gamma(d)$ is a totally geodesic submanifold of $M$ we have to show that
\[ \langle \nabla_X X , \eta(d) \rangle = 0 \, \, .\]

Let $F(s,t) := Exp(sX).\gamma(t)$ be the surface generated by the Killing vector field $X$ and the geodesic $\gamma(t)$.

Observe that:
\[ \begin{aligned} \langle \nabla_X X , \eta(t) \rangle &= \langle \frac{  \mathrm{D} }{\partial s} \frac{\partial}{\partial s}F(s,t){|}_{s=0}, \eta(t) \rangle  \\
\end{aligned} \]

So \[ \begin{aligned} \frac{\mathrm{d}}{\mathrm{d}t} \langle \nabla_X X , \eta(t) \rangle &= \langle \frac{ \mathrm{D} }{\partial t} \frac{ \mathrm{D} }{\partial s} \frac{\partial}{\partial s} F(s,t) |_{s=0}, \eta(t) \rangle  \\
&= \langle \frac{ \mathrm{D} }{\partial s} \frac{ \mathrm{D} }{\partial t} \frac{\partial}{\partial s} F(s,t) |_{s=0}, \eta(t) \rangle +
\langle R_{\gamma'(t),X}X, \eta(t) \rangle\\
&= 0 \, \, ,
\end{aligned} \]
\noindent
where the last equality follows from equations (\ref{dos-pajaros}) as explained above. Since $\langle \nabla_X X , \eta(0) \rangle = 0$ it follows
that $\langle \nabla_X X , \eta(t) \rangle \equiv 0$.

So $N$ is a totally geodesic submanifold of $M$ and this proves the proposition. $\Box$

\vspace{1cm}

Recall that the Karpelevich-Mostow theorem claims that a semisimple Lie subgroup $G \subset \mathrm{Iso}(M)$, where $M$ is a symmetric space of negative Ricci curvature, has a totally geodesic orbit $G.p \subset M$.

To prove the theorem notice that any such symmetric space $M$ can be equivariant and totally geodesically embedded in some $\mathcal{P} = SL(n,\mathbb{R})/SO(n)$. By the above proposition it is enough to show that a semisimple Lie subgroup $G \subset SL(n,\mathbb{R})$ has a totally geodesic orbit in $\mathcal{P}$. But this is a consequence of the so-called unitary Weyl's trick (for details see Proposition 2.2 of \cite{nostro}).

\section*{Karpelevich's proof}

After we wrote this corrigendum we found that the structure of the original proof by Karpelevich contained in \cite{Ka}
is the same as our proof. More precisely, Lemma 2 in \cite[page 402]{Ka} is the special case of Proposition \ref{out-in} in which
$G$ is assumed to be semisimple. It is clear that the hypothesis of $G$ being ``semisimple'' can not be omitted in Karpelevich's proof.
Thus, the idea of Karpelevich cannot be used to prove Proposition \ref{out-in}.

Then the proof of Theorem 1 \cite[page 403]{Ka} consists, as in our proof, in establishing the existence of a symmetric space $\mathcal{P}$ in which both semisimple Lie groups $G$ and $\mathrm{Iso}(M)$ have totally geodesic orbits. Notice that we took $\mathcal{P} = SL(n,\mathbb{R})/SO(n)$ whilst Karpelevich chose $\mathcal{P}$ as the symmetric space of non-compact type such that the Lie algebra of the isometry group $Iso(\mathcal{P})$ is the complexification of the Lie algebra of the Lie group $\mathrm{Iso}(M)$.

\section{Acknowledgements}

We thank Valeria Chiado' Piat for the help with the translation of \cite{Ka}.

{\small

\medskip
\noindent {\bf Authors '  Addresses:}\\
\noindent A. J. Di Scala, \\ Dipartimento di Matematica, \\
Politecnico di Torino, \\
Corso Duca degli Abruzzi 24, 10129 Torino, Italy \\
  email:     antonio.discala@polito.it \\

\noindent C. Olmos, \\ Fa.M.A.F.,\\ Universidad Nacional de C\'ordoba,
\\ Ciudad Universitaria,\\ 5000 Cordoba, Argentina\\
          email:          olmos@mate.uncor.edu\\

\end{document}